\documentclass[12pt]{amsart}
\usepackage{enumerate}

\theoremstyle{plain}
\newtheorem{theorem}{Theorem}
\newtheorem{lemma}{Lemma}

\newcommand{\Hs}{{\mathcal H}}

\newenvironment{enumeratei}{\begin{enumerate}[\upshape (i)]}{\end{enumerate}}

\begin{document}

\title{On capability of finite abelian groups}

\author{Zoran \v Suni\'c}

\address{Department of Mathematics, Texas A\&M University, College Station, TX 77843-3368, USA}

\email{sunic@math.tamu.edu}

\thanks{The work presented here is partially supported by NSF/DMS-0805932}


\begin{abstract}
Baer characterized capable finite abelian groups (a group is capable if it is isomorphic to the quotient of some group by its center) by a condition on the size of the factors in the invariant factor decomposition (the group must be noncyclic and the top two invariant factors must be equal). We provide a different characterization, given in terms of a condition on the lattice of subgroups. Namely, a finite abelian group $G$ is capable if and only if there exists a family $\{H_i\}$ of subgroups of $G$ with trivial intersection, such that the union generates $G$ and all the quotients $G/H_i$ have the same exponent. The condition that the family of subgroups generates $G$ may be replaced by the condition that the family covers $G$  and the condition that the quotients have the same exponent may be replaced by the condition that the quotients are isomorphic.  
\end{abstract}

\maketitle

A class of finite groups with a certain property played a crucial role in the construction of a large family of finitely generated torsion groups of intermediate growth~\cite{bartholdi-s:wpg} that generalize the examples of Grigorchuk~\cite{grigorchuk:gdegree}. Namely, Bartholdi and the author used finite groups $B$ that satisfy the following condition. There exists a family of normal subgroups $\{N_i\}_{i \in I}$ of $B$ such that 
\begin{enumeratei}
\item $\bigcap_{i \in I} N_i = 1$,

\item $\bigcup_{i \in I} N_i = B$,

\item $B/N_i \cong B/N_j$, for $i,j \in I$. 
\end{enumeratei}
In his dissertation~\cite{sunik:phd} the author conjectured that the class of finite abelian groups that have this property is precisely the class of noncyclic abelian groups for which the top two factors in the invariant factor decomposition are equal. We prove this conjecture here. By an earlier result of Baer~\cite{baer:capable-abelian}, it follows that this class of finite abelian groups is precisely the class of capable finite abelian groups. This gives a complete description of the finite abelian groups $B$ that satisfy the condition on the lattice of subgroups given above. The case of nonabelian groups seems much harder and not much is known at the moment. Note that, because of the intersection condition, the group $B$ embeds as a subdirect product in the product $\prod_{i\in I} B/N_i$. Therefore, $B$ is nonabelian if and only if the quotients $B/N_i$ are nonabelian. The simplest example of a nonabelian group $B$ with the above property of which the author is aware uses the symmetric group $S_3$ for the quotients and has size $4 \cdot 3^6=2916$. It was communicated to Bartholdi and the author by D.~Holt (see~\cite{bartholdi-s:wpg} for specifics on this nonabelian example). A necessary condition on the structure of $B$ is given by Brodie, Chamberlain, and Kappe~\cite{brodie-c-k:covers}. Namely, they show that a group can be covered by a family of proper normal subgroups if and only if it has a noncyclic abelian quotient. 

Note that the condition on the lattice of subgroups given above implies the condition (c) in the following theorem and is implied by the condition (d). 

\begin{theorem} \label{t:baerconditions}
Let $G \cong C_{n_1} \times \dots \times C_{n_k}$ be a finite nontrivial abelian group, where $C_{n_i}$, $i=1,\dots,k$, denotes the cyclic group of order $n_i$ and $1< n_1 \mid n_2 \mid \dots \mid n_k$.

The following conditions are equivalent.
\begin{itemize}

\item[\textup{(a)}]
There exists a group $L$ such that $G \cong L/Z(L)$, where $Z(L)$ denotes the center of $L$.

\item[\textup{(b)}] $k \geq 2$ and $n_{k-1}=n_k$.

\item[\textup{(c)}]
There exists a family of subgroups $\{H_i\}_{i \in I}$ of $G$ such that
\begin{enumeratei}
\item $\bigcap_{i \in I} H_i = 1$,

\item $\langle \bigcup_{i \in I} H_i \rangle = G$,

\item all groups in the family $\{G/H_i\}_{i \in I}$ have the same exponent. 

\end{enumeratei}

\item[\textup{(d)}]
There exists a family of subgroups $\{H_i\}_{i \in I}$ of $G$ such that
\begin{enumeratei}
\item $\bigcap_{i \in I} H_i = 1$,

\item $\bigcup_{i \in I} H_i = G$,

\item $G/H_i \cong G/H_j$, for $i,j \in I$.

\item $H_i \cong H_j$, for $i,j \in I$.
\end{enumeratei}

\end{itemize}
\end{theorem}

Before we prove the main result, we provide a statement of a simple lemma that is certainly known. Nevertheless, for completeness, we provide an elementary proof of the lemma after the proof of Theorem~\ref{t:baerconditions}. 

\begin{lemma}\label{lemma}
Let $A=\langle a \rangle \cong B = \langle b \rangle \cong C_n$, where $n>1$, and let $C=AB$ be the direct product of $A$ and $B$ (written internally). Let $X$ be the set of elements of $C$ of order $n$. 
\begin{itemize}
\item[\textup{(a)}] $a^ib^j \in X$ if and only if the greatest common divisor of $i$, $j$, and $n$ is 1. 

\item[\textup{(b)}] For every $x \in X$ there exists $y \in X$ such that $C=\langle x \rangle \langle y \rangle$ as an internal direct product. 

\item[\textup{(c)}] $C = \bigcup_{x \in X} \langle x \rangle$.  
\end{itemize}  
\end{lemma} 

\begin{proof}[Proof of Theorem~\ref{t:baerconditions}]
Let $\{G_i\}_{i=1}^k$ be a family of subgroups of $G$ such that,
$\langle a_i \rangle = G_i \cong C_{n_i}$, for $i=1,\dots,k$ and
$G=G_1G_2\dots G_k$. \smallskip

(a) is equivalent to (b). This is proved by Baer~\cite{baer:capable-abelian}.  \smallskip

(d) implies (c). Clear. \smallskip

(c) implies (b). Assume that the family $\{H_i\}_{i \in I}$ of subgroups of $G$ satisfies the intersection condition (c)(i), the generating condition (c)(ii), and the quotient exponent condition (c)(iii). 

By way of contradiction, assume that $G$ is cyclic. Then all of its subgroups and quotients are cyclic and, for every subgroup $H$ of $G$,  the exponent of $G/H$ is equal to $|G|/|H|$. Since, for every divisor of $|G|$, there exists a unique subgroup of that particular order, the quotient exponent condition implies that the family $\{H_i\}_{i \in I}$ consists of a single group, which, by itself, cannot satisfy both the intersection condition and the generating condition in a nontrivial group. Therefore $G$ is not cyclic, i.e., $k \geq 2$. 

By way of contradiction, assume that $n_k/n_{k-1} = m >1$. 

Let $p$ be a prime dividing $m$ and $A_p$ be the unique subgroup of $G_k$ of order $p^t$, where $p^t$ is the highest power of $p$ dividing $m$. 

Let $H$ be a subgroup of $G$. We claim that $G/H$ has index dividing $\frac{n_k}{p^t}=n_{k-1}\cdot \frac{m}{p^t}$ if and only if $H$ contains $A_p$. 

Assume that $H$ contains $A_p$. For any element $g$ in $G$, the element $g^{n_{k-1}}$ belongs to $G_k$. Therefore $g^{n_k/p^t}$ is an element of $G_k$ of order dividing $p^t$, which means that $g^{n/p^t}$ belongs to $A_p \leq H$. Thus $g^{n_k/p^t}=1$ in $G/H$, showing that the exponent of $G/H$ divides $n_k/p^t$. 

Conversely, let the exponent of $G/H$ divide $n_k/p^t$. Let $p^T$ be the highest power of $p$ dividing $n_k$ and let $g$ be an element of $G_k$ of order $p^T$. The element $g^{n/p^t}$ must be equal to 1 in $G/H$, which means that $g^{n_k/p^t}$ belongs to $H$. On the other hand, $g^{n_k/p^t}$ is an element of $G_k$ of order exactly $p^t$, which means that it generates $A_p$. Therefore $H$ contains $A_p$. 

Since the exponent of $G$ is $n_k$, the generating condition for the family $\{H_i\}_{i \in I}$ implies that at least one group in this family must have exponent divisible by $p^T$ (since the exponent of the group generated by the family  $\{H_i\}_{i \in I}$ is the least common multiple of the exponents of the groups $H_i$, $i \in I$). Let $H_j$ be such a group. Since $H_j$ is abelian, it must contain an element $h$ of order $p^T$.  The element $h^{p^{T-t}}$ is then an element of $G_k$ of order $p^t$, which means that it generates $A_p$. Therefore $H_j$ contains $A_p$ and the exponent of $G/H_j$ divides $n_k/p^t$. 

By the last remark and the quotient exponent condition, $G/H_i$ divides $n/p^t$, for all $i \in I$, which implies that each member of the family $\{H_i\}_{i \in I}$ contains $A_p$, violating the intersection condition.

Therefore, $k \geq 2$ and $n_k=n_{k-1}$.  \smallskip

(b) implies (d). Let $H_i=\langle a_1,\dots,a_{i-1},a_ia_k^{n_k/n_i},a_{i+1},\dots,a_{k-1} \rangle$, $i=1,\dots,k-1$, and $H_k=\langle a_1,\dots,a_{k-1} \rangle$.

Clearly $H_k \cong G_1\dots G_{k-1}$. Further, for $i=1,\dots,k-1$, the order of the element $a_ia_k^{n_k/n_i}$ is exactly $n_i$. Since $\langle a_ia_k^{n_k/n_i} \rangle \cap \langle a_1,\dots,a_{i-1},a_{i+1},\dots,a_{k-1} \rangle = 1$ we have that, for $i=1,\dots,k$,  
\[ H_i \cong G_1G_2 \dots G_{k-1}. \]

Clearly, $G/H_k \cong G_k$. We claim that $G/H_i \cong G_k$, for $i=1,\dots,k-1$. We have $G/H_i \cong G_iG_k/\langle a_ia_k^{n_k/n_i}\rangle$. Since $a_ia_k^{n_k/n_i},a_k \in \langle a_ia_k^{n_k/n_i} \rangle G_k$ we conclude that $G_iG_k = \langle a_ia_k^{n_k/n_i} \rangle G_k$. Therefore $G/H_i \cong \langle a_ia_k^{n_k/n_i} \rangle G_k/\langle a_ia_k^{n_k/n_i} \rangle$. We claim that $G_k \cap \langle a_ia_k^{n_k/n_i} \rangle= 1$. Indeed, if $g = (a_ia_k^{n_k/n_i})^m \in G_k$, then $n_i$ (the order of $a_i$) must divide $m$, which implies that $g = a_i^m a_k^{mn_k/n_i}=1$. Thus the intersection $\langle a_ia_k^{n_k/n_i} \rangle \cap G_k$ is trivial and we have $\langle a_ia_k^{n_k/n_i} \rangle G_k/\langle a_ia_k^{n_k/n_i} \rangle \cong G_k$. Thus, for $i=1,\dots,k$, 
\[ G/H_i \cong G_k.  \]

Let $g=a_1^{m_1} \dots a_k^{m_k}$, $0 \leq m_i < n_i$ for $i=1,\dots,k$, be an element in the intersection $\cap_{i=1}^k H_i$. Since $g \in H_k=\langle a_1,\dots,a_{k-1} \rangle$, we must have $m_k=0$. For $i=1,\dots,k-1$, since $g \in H_i=\langle a_1,\dots,a_{i-1},a_ia_k^{n_k/n_i},a_{i+1},\dots,a_{k-1} \rangle$ and $m_k=0$, the exponent $m_i$ must be divisible by $n_i$ (the order of the element $a_k^{n_k/n_i}$), which then implies that $m_i=0$. Thus, 
\[ \bigcap_{i-1}^k H_i = 1. \]

Since $a_1,\dots,a_{k-1} \in H_k$ and $a_{k-1}a_k \in H_{k-1}$, we have $\langle \bigcup_{i-1}^k H_i \rangle = G$. 

So far, we constructed a family of subgroups $\Hs = \{H_i\}_{i=1}^k$ of $G$ that are all isomorphic to $G_1\dots G_{k-1}$, with quotients isomorphic to $G_k$, with trivial intersection, and with union that generates $G$. However, we need a family of groups that cover $G$. In order to accomplish this we will enlarge the family $\Hs$ by more subgroups $H$ of $G$ such that $H \cong G_1\dots G_{k-1}$ and $G/H \cong G_k$. 

Denote $n_{k-1}=n_k=n$, $G_{k-1}=A=\langle a\rangle$ and $G_k = B = \langle b \rangle$. Let $X$ be the set of elements of order $n$ in $AB$. For $x$ in $X$, let $H_x=G_1\dots G_{k-2}\langle x \rangle$. Clearly
\[ H_x \cong G_1 \dots G_{k-1}. \]

By part (b) of Lemma~\ref{lemma}, 
\[ G/H_x \cong G_1 \dots G_{k-2}AB/G_1\dots G_{k-2}\langle x \rangle \cong AB/\langle x \rangle \cong C_n \cong G_k. \]

By part (c) of Lemma~\ref{lemma}, 
\[ 
 \bigcup_{x \in X} H_x = \bigcup_{x \in X} G_1 \dots G_{k-2}\langle x \rangle = 
  G_1 \dots G_{k-2} \bigcup_{x \in X} \langle x \rangle = G_1 \dots G_{k-2}AB = G.
\]

Therefore, the union of the families $\{H_i\}_{i=1}^k$ and $\{H_x\}_{x \in X}$ satisfies the condition (d).  
\end{proof}

\begin{proof}[Proof of Lemma~\ref{lemma}]
In what follows, given an element $a^ib^j$ in $C$, let $d$ be the greatest common divisor of $i$ and $j$ and let $i'$ and $j'$ be integers such that $ii'+jj'=d$.   

(a) If $(a^ib^j)^m = 1$, then $mi \equiv mj \equiv 0 \pmod{n}$, implying $md \equiv m(ii'+jj') \equiv 0 \pmod{n}$. If $d$ is relatively prime to $n$, then we must have $m \equiv 0 \pmod {n}$. In other words, if $i$, $j$, and $n$ are relatively prime, the order of $a^ib^j$ is $n$. On the other hand, if $i$, $j$, and $n$ have a common divisor $d'>1$ then $(a^ib^j)^{n/d'}=1$, implying that the order of $a^ib^j$ is smaller than $n$.  

(b) Let $x=a^ib^j \in X$. Since $x \in X$, $d$ is relatively prime to $n$. Consider the element $y=a^{j'}b^{-i'}$. We have $x^{i'}y^{j} = a^{ii'}b^{ji'}a^{jj'}b^{-ji'}= a^d$. Similarly, $x^{j'}y^{-i}=b^d$. Since $d$ is relatively prime to $n$, $a^d$ and $b^d$ generate $C$, implying that $x$ and $y$ also generate $C$. Since the order of $y$ cannot be larger than $n$, it must be equal to $n$ and $C$ is internal direct product of $\langle x \rangle$ and $\langle y \rangle$.  

(c) Let $c=a^ib^j$ be a nontrivial element of $C$. Then $c=a^i b^j = (a^{i''} b^{j''})^d$, where $i''$ and $j''$ are relatively prime. Therefore $c \in \langle c'' \rangle$, where $c''=a^{i''}b^{j''}$ has order $n$. We conclude that $C=\bigcup_{x \in X} \langle x \rangle$.
\end{proof}


\begin{thebibliography}{BCK88}

\bibitem[Bae38]{baer:capable-abelian}
Reinhold Baer.
\newblock Groups with preassigned central and central quotient group.
\newblock {\em Trans. Amer. Math. Soc.}, 44(3):387--412, 1938.

\bibitem[BCK88]{brodie-c-k:covers}
M.~A. Brodie, R.~F. Chamberlain, and L.-C. Kappe.
\newblock Finite coverings by normal subgroups.
\newblock {\em Proc. Amer. Math. Soc.}, 104(3):669--674, 1988.

\bibitem[B{\v{S}}01]{bartholdi-s:wpg}
Laurent Bartholdi and Zoran {\v{S}}uni{\'k}.
\newblock On the word and period growth of some groups of tree automorphisms.
\newblock {\em Comm. Algebra}, 29(11):4923--4964, 2001.

\bibitem[Gri84]{grigorchuk:gdegree}
R.~I. Grigorchuk.
\newblock Degrees of growth of finitely generated groups and the theory of
  invariant means.
\newblock {\em Izv. Akad. Nauk SSSR Ser. Mat.}, 48(5):939--985, 1984.

\bibitem[{\v{S}}un00]{sunik:phd}
Zoran {\v{S}}uni{\'k}.
\newblock {\em On a class of periodic spinal groups of intermediate growth}.
\newblock PhD thesis, State University of New York at Binghamton, 2000.

\end{thebibliography}

\def\cprime{$'$}

\end{document}